\documentclass[12pt]{article}
\usepackage{amsmath,amssymb,theorem}
\usepackage{diagbox}
\usepackage{multirow}
\usepackage{graphicx}
\usepackage{subfigure}
\usepackage{psfrag}
\usepackage{color}
\usepackage{enumerate}
\usepackage{epstopdf}
\usepackage[T1]{fontenc}
\usepackage[colorlinks,linkcolor=blue]{hyperref}
\newtheorem{thm}{Theorem}[section]
\newtheorem{conj}[thm]{Conjecture}
\newtheorem{cor}[thm]{Corollary}
\newtheorem{prop}[thm]{Proposition}
\newtheorem{que}[thm]{Question}
\theorembodyfont{\rmfamily}
\newtheorem{defn}[thm]{Definition}
\newtheorem{obser}[thm]{Observation}
\def\pf{\bigskip\noindent {\bf Proof.}~~}

\def\dfn#1{{\sl #1}}

\def\less{\backslash}

\def\pf{\bigskip\noindent {\bf{Proof.}}~~}

\def\fattextindent#1{\indent\indent\llap{#1\enspace}\ignorespaces}

\def\myitemitem{\par\hangindent\parindent\fattextindent}

\title{A survey on star edge-coloring of graphs}
\author{Hui Lei$^1$, Yongtang Shi$^2$\\
{$^1$ School of Statistics and Data Science, LPMC and KLMDASR}\\
{ Nankai University, Tianjin 300071, China}\\
{ $^2$ Center for Combinatorics and LPMC}\\
{ Nankai University, Tianjin 300071,  China}\\
 {hlei@nankai.edu.cn; shi@nankai.edu.cn}\\
}
\date{}

\begin{document}
\maketitle

\begin{abstract}
The \dfn{star chromatic index}  of a multigraph $G$, denoted
$\chi'_{st}(G)$,  is the minimum number of colors needed to properly
color the edges of $G$ such that no path or cycle of length four is
bicolored.  We survey the results of determining
the star chromatic index, present the  interesting proofs and techniques, and collect many open problems
and conjectures.  \\

\end{abstract}

\section{Introduction}

\baselineskip 17pt

All multigraphs  in this paper are finite and loopless; and all graphs  are finite and without loops or multiple edges.
Given a multigraph $G$, let $c: E(G)\rightarrow [k]$ be a proper
edge-coloring of $G$, where $k\ge1$ is an integer and $[k]:=\{1,2, \dots, k\}$. We say that
$c$ is a  \dfn{star $k$-edge-coloring} of $G$ if no path or cycle of
length four in $G$ is bicolored under the coloring $c$; and  $G$ is
\dfn{star $k$-edge-colorable} if $G$ admits a star
$k$-edge-coloring. The \dfn{star chromatic index}  of $G$, denoted $\chi'_{st}(G)$,  is the smallest integer $k$ such that $G$ is
star $k$-edge-colorable.  The  chromatic index and chromatic number of $G$ are denoted by  $\chi'(G)$ and $\chi(G)$.
As pointed out in  \cite{DMS2013}, the definition of
star edge-coloring of a graph $G$ is equivalent to the star
vertex-coloring of its line graph $L(G)$.  Star edge-coloring of a
graph was initiated  by Liu and Deng \cite{DL2008}, motivated by the
vertex version (see \cite{ACKKR2004, BCMRW2009, CRW2013, CWZ2019, KKT2009,
NM2003}).   Given a multigraph $G$, we use  $|G|$ to denote the number of vertices,  $e(G)$ the number of edges, $\delta(G)$ the minimum degree, and $\Delta(G)$ the maximum degree of $G$, respectively.  
We use $K_n$, $P_n$ and $C_n$ to denote the complete graph, the path and the cycle on $n$ vertices, respectively. The \dfn{maximum average degree} of  a multigraph $G$, denoted  $\text{mad}(G)$, is defined as the maximum  of  $2
e(H)/|H|$ taken over all the non-empty subgraphs $H$ of $G$.  The  {\em girth} of a graph with a cycle is the length of its shortest cycle. A graph with no cycle has infinite girth. The following first upper bound is a result of  Liu and Deng \cite{DL2008}.  \medskip

\begin{thm}[\cite{DL2008}]
Every graph  $G$ with  $\Delta(G)\geq7$ satisfies
$\chi'_{st}(G)\leq \lceil16(\Delta(G)-1)^\frac3 2\rceil.$
\end{thm}

Theorem~\ref{Kn}  below is a result of Dvo\v{r}\'ak, Mohar, and \v{S}\'amal~\cite{DMS2013},   which give an   upper and a lower bounds for complete graphs.

\begin{thm} [\cite{DMS2013}]\label{Kn}  The star chromatic index of   the complete graph $K_n$ satisfies

$$2n(1+o(1))\leq \chi'_{st}(K_n)\leq n\, \frac{2^{2\sqrt{2}(1+o(1))\sqrt{\log n}}}{(\log n)^{1/4}}.$$
In particular, for every $\epsilon>0$, there exists a constant $c$ such that  $\chi'_{st}(K_n)\le cn^{1+\epsilon}$ for every integer $n\ge1$.
\end{thm}

They proved the upper bound using a nontrivial result about sets without arithmetic progressions, and up till now, it is still the best known. For the lower bound, they used an elegant double counting approach. The authors of \cite{BLM2016} observed a improvement in their proof and obtained the bound $\chi'_{st}(K_n)\geq3n(n-1)/(n+4)$ (see \cite{M2013} for a proof).

Currently the best upper bound in terms of the maximum degree for general  graphs is also derived in~\cite{DMS2013} using Theorem~\ref{Kn}.
\begin{thm}[\cite{DMS2013}]
For  any  graph $G$ with maximum degree $\Delta$, $\chi'_{st}(G)\leq \Delta\cdot 2^{O(1)\sqrt{\log \Delta}}$.
\end{thm}

The true order of magnitude of  $\chi'_{st}(K_n)$ is still unknown. Dvo\v{r}\'ak, Mohar, and \v{S}\'amal~\cite{DMS2013} made the following question.
\begin{que}[\cite{DMS2013}]
What is the true order of magnitude of $\chi'_{st}(K_n)$? Is $\chi'_{st}(K_n)=O(n)$?
\end{que}

Dvo\v{r}\'ak, Mohar, and \v{S}\'amal~\cite{DMS2013} also mentioned that, if
convenient, one may start with other special graphs in place of $K_n$, in particular, from Observation \ref{knn} it follows that if $\chi'_{st}(K_{n,n})$ is $O(n)$(or $n(\log n)^{O(1)}$, $n^{1+o(1)}$, respectively), then $\chi'_{st}(K_n)$ is $O(n\log n)$(or $n(\log n)^{O(1)}$, $n^{1+o(1)}$, respectively).
\begin{obser}[\cite{DMS2013}]\label{knn}
$\chi'_{st}(K_{n})\leq\sum\limits_{i=1}^{\lceil\log_2n\rceil}2^{i-1}\chi'_{st}(K_{\lceil n/2^i\rceil,\lceil n/2^i\rceil})$.
\end{obser}

The paper is organized as follows. The first Section \ref{Complexity} contains the results on
algorithmic aspects of the star edge-coloring.  In Sections \ref{Subcubicgraphs}, \ref{Bipartitegraphs}, \ref{Planargraphs} we  present the results of star chromatic index on  subcubic graphs, bipartite graphs, and planar graphs, respectively. Star edge-coloring is naturally generalized to the list version and  it was pointed out in \cite{DMS2013}: It would be
interesting to understand the list version of star edge-coloring. The results of list star edge-coloring are presented in  Section \ref{Liststaredgecoloring}.
Finally in Section \ref{Anotherniceclasses} we summarize the results of  the star chromatic index of another family of graphs such as  Halin graphs,  generalized Petersen graphs,  products of graphs and so on.

\section{Complexity}\label{Complexity}

It  is well-known~\cite{Vizing}  that the chromatic index of a graph with maximum degree $\Delta$ is   either $\Delta$ or $\Delta+1$. However,  it is NP-complete~\cite{Holyer} to determine whether the chromatic index of an
arbitrary graph with maximum degree $\Delta$ is $\Delta$ or $\Delta+1$.  The problem remains NP-complete even for cubic graphs (the degree of all vertices is 3).   Lei, Shi, and Song \cite{LSS2017}  proved the following result.

\begin{thm}[\cite{LSS2017}]\label{NP}
It is  NP-complete to determine whether $\chi'_{st}(G)\le3$ for an arbitrary graph $G$.
\end{thm}
\pf First let us denote by SEC the  problem stated in Theorem~\ref{NP}, and   we denote by  3EC the following well-known NP-complete problem of Holyer~\cite{Holyer}: \medskip

Given a  cubic graph $G$, is $G$ $3$-edge-colorable? \medskip

 Clearly, SEC is  in the class   NP.  We shall reduce 3EC  to SEC. \medskip

\begin{figure}[htbp]
\begin{center}
\scalebox{1.2}[1.2]{\includegraphics{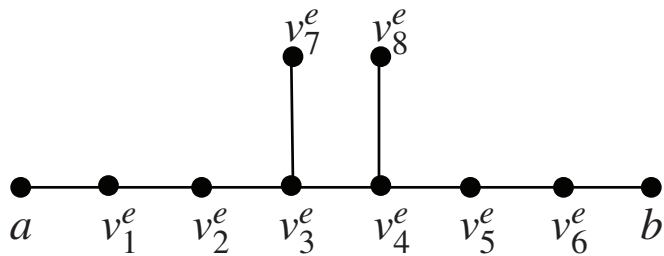}}
\caption{Graph $H_{ab}$.}\label{fig1}
\end{center}
\end{figure}

Let  $H$ be an instance of 3EC. We construct a   graph $G$ from $H$ by replacing   each edge $e=uw\in E(H)$  with a copy of  graph $H_{ab}$, identifying $u$ with $a$ and $w$ with $b$, where $H_{ab}$ is depicted in Figure \ref{fig1}.
The size of $G$ is clearly polynomial in the size of $H$, and $\Delta(G)=3$.

It suffices to show that   $\chi' (H)\le3$ if and only if $\chi'_{st}(G)\le3$.  Assume that $\chi'(H)\le3$.  Let $c:E(H)\rightarrow \{1,2,3\}$
be a proper $3$-edge-coloring of $H$. Let  $c^*$  be an edge coloring of $G$ obtained from $c$  as follows:
for each edge $e=uw\in E(H)$, let  $c^*(av^e_1)=c^*(v^e_3v^e_4)=c^*(v^e_6b)=c(uw)$,
$c^*(v^e_1v^e_2)=c^*(v^e_3v^e_7)=c^*(v^e_4v^e_8)=c^*(v^e_5v^e_6)=c(uw)+1$,
and $c^*(v^e_2v^e_3)=c^*(v^e_4v^e_5)=c(uw)+2$, where all colors here and henceforth are done modulo $3$.
Notice  that $c^*$ is a proper $3$-edge-coloring of $G$. Furthermore,  it can be easily checked  that
$G$ has no bicolored path or  cycle of length four under the coloring $c^*$. Thus  $c^*$ is a star $3$-edge-coloring of $G$ and so  $\chi'_{st}(G)\le3$.\medskip

Conversely, assume that  $\chi'_{st}(G)\le3$. Let  $c^*:E(G)\rightarrow \{1,2,3\}$
be a star $3$-edge-coloring of $G$.  Let $c$ be an edge-coloring of $H$ obtained from $c^*$ by letting $c(e)=c^*(av^e_1)$ for any $e=uw\in E(H)$.  Clearly, $c$ is a proper $3$-edge-coloring of $H$ if  for any edge  $e=uw$ in $G$, $c^*(av^e_1)=c^*(v^e_6b)$. We prove this next.   Let $e=uw$ be an edge of $H$. We consider the following two cases. \medskip

\noindent {\bf Case 1:} $c^*(v^e_3v^e_7)=c^*(v^e_4v^e_8)$.\medskip

In this case, let  $c^*(v^e_3v^e_7)=c^*(v^e_4v^e_8)=\alpha$, where $\alpha\in\{1,2,3\}$. We may further assume that  $c^*(v^e_3v^e_4)=\beta$ and $c^*(v^e_2v^e_3)=c^*(v^e_4v^e_5)=\gamma$, where $\{\beta, \gamma\}=\{1,2,3\}\less \alpha$. This is possible because   $d_{G}(v^e_3)=d_{G}(v^e_4)=3$  and $c^*$ is a proper $3$-edge-coloring of $G$.  Since $c^*$ is a star edge-coloring of $G$, we see that $c^*(v^e_1v^e_2)=c^*(v^e_5v^e_6)=\alpha$ and so $c^*(av^e_1)=c^*(v^e_6b)=\beta$. \medskip

\noindent {\bf Case 2:} $c^*(v^e_3v^e_7)\ne c^*(v^e_4v^e_8)$.\medskip

In this case, let  $c^*(v^e_3v^e_7)=\alpha$, $c^*(v^e_4v^e_8)=\beta$, $c^*(v^e_3v^e_4)=\gamma$, where $\{\alpha, \beta, \gamma\}=\{1,2,3\}$. This is possible because $\alpha\ne \beta$ by assumption.
Since $c^*$ is a proper edge-coloring of $G$, we see that   $c^*(v^e_2v^e_3)=\beta$  and $c^*(v^e_4v^e_5)=\alpha$. One can easily check now that $c^*(v^e_1v^e_2)=\alpha$ and $c^*(v^e_5v^e_6)=\beta$,  and so $c^*(av^e_1)=c^*(v^e_6b)=\gamma$, because $c^*$ is a star edge-coloring of $G$. \medskip

In both cases we see that $c^*(av^e_1)=c^*(v^e_6b)$. Therefore  $c$ is  a proper $3$-edge-coloring of $H$ and so $\chi'(H)\le3$.  This completes the proof of Theorem~\ref{NP}.
\hfill\vrule height3pt width6pt depth2pt\\

There are few results for the computational complexity of the star edge-coloring problem. In \cite{ORD2018}, Omoomi, Roshanbin, and Dastjerdi  presented a polynomial time algorithm that finds an
optimum star edge-coloring for every tree.
\begin{thm}[\cite{ORD2018}]
There is a polynomial time algorithm for computing the star chromatic index of
every tree and presenting an optimum star edge-coloring of it.
\end{thm}
\section{Subcubic graphs}\label{Subcubicgraphs}
 A multigraph $G$ is \dfn{subcubic} if all its vertices have degree less than or equal to
three.  Dvo\v{r}\'ak, Mohar, and \v{S}\'amal~\cite{DMS2013} considered the star chromatic index of subcubic multigraphs. To state their result, we need to introduce one notation. A graph $G$ \dfn{covers} a graph $H$ if there is a mapping $f: V(G)\rightarrow V(H)$ such that for any $uv\in E(G)$, $f(u)f(v)\in E(H)$,
 and for any $u\in V(G)$, $f$ is a bijection between $N_G(u)$ and $N_{H}(f(u))$.  They proved the following.

\begin{thm} [\cite{DMS2013}]\label{s=7} Let $G$ be a multigraph.
\begin{enumerate}[(a)]
\item  If $G$ is  subcubic,  then $\chi'_{st}(G)\le7$.\vspace{-8pt}

\item   If $G$ is  cubic and has no multiple edges, then $\chi'_{st}(G)\ge4$ and the equality holds if and only if $G$ covers the graph of $3$-cube.
\end{enumerate}
\end{thm}

As observed in~\cite{DMS2013},  $K_{3,3}$ and  the Heawood graph are star $6$-edge-colorable. No subcubic multigraphs with star chromatic index seven  are known. Dvo\v{r}\'ak, Mohar, and \v{S}\'amal~\cite{DMS2013}   proposed the following conjecture.

\begin{conj} [\cite{DMS2013}]\label{cubic}
Let $G$ be a subcubic multigraph. Then $\chi'_{st}(G)\leq 6$.
\end{conj}

  It was  shown  that
every subcubic outerplanar graph is star $5$-edge-colorable in~\cite{BLM2016} and every cubic Halin graph is star $6$-edge-colorable in~\cite{CGR2019,HLW2020}.
 Lei, Shi, and Song~\cite{LSS2017} proved that
\begin{thm}[\cite{LSS2017}]\label{thm*}
Let $G$ be a subcubic multigraph.\medskip

\myitemitem {(a)}   If $mad(G)<2$, then $\chi'_{s}(G)\leq 4$ and the bound is tight. \medskip

\myitemitem {(b)}  If $mad(G)<24/11$, then $\chi'_{s}(G)\leq 5$.\medskip

\myitemitem {(c)} If $mad(G)<5/2$, then $\chi'_{s}(G)\leq 6$.
\end{thm}

 \begin{figure}[htbp]
\begin{center}
\includegraphics[scale=0.4]{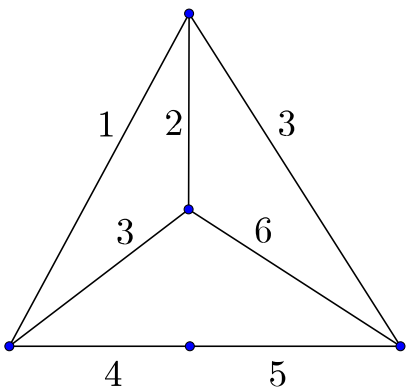}
\caption{ A graph  with maximum average degree $14/5$ and star chromatic index  $6$.}\label{prism}
\end{center}
\end{figure}

Later,  Lei, Shi, Song, and Wang~\cite{LSSW2018} improved Theorem~\ref{thm*}(b) by showing the following  result.

\begin{thm}[\cite{LSSW2018}]\label{mainthm}
Let $G$ be a subcubic multigraph with $mad(G)<12/5$. Then $\chi'_{st}(G)\leq 5$.
\end{thm}

  We don't know if the bound $12/5$  in Theorem~\ref{mainthm} is best possible.  The graph depicted  in Figure~\ref{prism}  has  maximum average degree $14/5$  but is not star $5$-edge-colorable.  \medskip

In \cite{WWW2019}, Wang, Wang, and Wang focused on the star edge-coloring of graphs with maximum
degree four and proved the following result.
\begin{thm}[\cite{WWW2019}]
Every graph $G$ with $\Delta(G)=4$ satisfies $\chi'_{st}(G) \le14.$
\end{thm}

\section{Bipartite graphs}\label{Bipartitegraphs}
In this section we consider the star edge-coloring of bipartite graphs. We first consider complete
bipartite graphs. Let $K_{m,n}$ denote the complete bipartite graph in which the orders of its bipartition sets are $m$ and $n$, where $m\leq n$. Trivially $\chi'_{st}(K_{1,n})=n$.

Dvorak,  Mohar,  and \v{S}\'amal \cite{DMS2013} obtained Observation \ref{obserknn}, thus proved an asymptotically  bound on $\chi'_{st}(K_{n,n})$: it follows from Theorem \ref{Kn} that for every $\varepsilon>0$, there is a constant $C>0$, such that for every $n\geq1$, $\chi'_{st}(K_{n,n})\leq Cn^{1+\varepsilon}$.
\begin{obser}[\cite{DMS2013}]\label{obserknn}
$\chi'_{st}(K_{n,n})\leq\chi'_{st}(K_n)+n$.
\end{obser}

Wang, Wang, and Wang \cite{WWW2019} proved that $\chi'_{st}(K_{3,4})=7$, and it is known that $\chi'_{st}(K_{3,3})=6$ in \cite{DMS2013}.
Casselgren, Granholm, and Raspaud \cite{CGR2019} obtained the following results.
\begin{thm}[\cite{CGR2019}]\label{2345}
Let $m, n$ be positive integers. The following holds.
\begin{enumerate}[(a)]
\item $\chi'_{st}(K_{2,n})=2n-\lfloor\frac{n}{2}\rfloor$.
\item $\chi'_{st}(K_{3,n})=3\lceil\frac{n}{2}\rceil$ for $n\neq4$.
\item $\frac{5n}{3}\leq\chi'_{st}(K_{4,n})\leq20\lceil\frac{n}{12}\rceil$.
\item $\chi'_{st}(K_{m,n})\leq15\lceil\frac{n}{8}\rceil\lceil\frac{m}{8}\rceil.$
\end{enumerate}
\end{thm}

Some exact values of $\chi'_{st}(K_{m,n})$ are given by Casselgren, Granholm, and Raspaud \cite{CGR2019} in the following Table \ref{tab:1}.


Moreover, using the idea in the proof of Theorem \ref{2345}(c), one can prove lower bounds on the star
chromatic index for further families of complete bipartite graphs. Let us here just list a few cases
corresponding to the values in the Table \ref{tab:1}(see \cite{CGR2019}).

\begin{itemize}
  \item $\chi'_{st}(K_{5,n})\geq\frac{5n}{3}$.
  \item $\chi'_{st}(K_{6,n})\geq\frac{7n}{4}$.
  \item $\chi'_{st}(K_{7,n})\geq\frac{7n}{4}$.
  \item $\chi'_{st}(K_{8,n})\geq\frac{9n}{5}$.
\end{itemize}

\begin{table}[htbp]
\begin{center}
\begin{tabular}{|c|c|c|c|c|c|c|c|c|c|}
  \hline
 $n$ &4& 5 &6 &7& 8& 9& 10& 11&12\\\hline
   $\chi'_{st}(K_{4,n})$ & 7&  10 & 11 &  13 &14 & 16 & 17& 20&20\\\hline
  $\chi'_{st}(K_{5,n})$ & &  11 & 12 &  14 &15 & 17 & 18& 20&\\\hline
  $\chi'_{st}(K_{6,n})$ & &  & 13 &  14 &15 & & & &\\\hline
  $\chi'_{st}(K_{7,n})$ & &  &  &  14 &15 & & & &\\\hline
  $\chi'_{st}(K_{8,n})$ & &  &  &   &15 & & & &\\\hline
\end{tabular}
\end{center}
\caption{The values of $\chi'_{st}(K_{m,n})$ for $m\in\{4,5,6,7,8\}$ and $n\leq12$.}\label{tab:1}
\end{table}
Next, we turn to general bipartite graphs with restrictions on the vertex degrees.
In the following we use the notation $G=(X,Y;E)$ for a bipartite graph $G$ with parts $X$ and
$Y$ and edge set $E=E(G)$. We denote by $\Delta(X)$ and $\Delta(Y)$ the maximum degrees of the vertices
in the parts $X$ and $Y$, respectively. A bipartite graph $G=(X,Y;E)$ where all vertices in $X$ have
degree $r$ and all vertices in $Y$ have degree $d$ is called $(r,d)$-biregular.

If the vertices in one part of a bipartite graph $G$ have maximum degree 1, then trivially $\chi'_{st}(G)=\Delta(G)$. For the case when the vertices in one of the parts have maximum degree two, Casselgren, Granholm, and Raspaud \cite{CGR2019} proved the following.
\begin{thm}[\cite{CGR2019}]\label{X2}
Let $G=(X,Y;E)$ be a bipartite graph with $\Delta(X)=2$ and $k\geq1$. Then
\begin{itemize}
  \item $\chi'_{st}(G)\leq3k$ if $\Delta(Y)=2k$.
  \item $\chi'_{st}(G)\leq3k+2$ if $\Delta(Y)=2k+1$.
\end{itemize}
\end{thm}
Note that the upper bounds in Theorem \ref{X2} are sharp. This can be obtained from Theorem \ref{2345}(a).

By Theorem \ref{X2}, one can get the following corollary for $(2,d)$-biregular graphs.
\begin{cor}[\cite{CGR2019}]
If $G$ is a $(2,d)$-biregular graph with $d\geq3$, then $\chi'_{st}(G)\leq\frac{3d+1}{2}$.
\end{cor}
Melinder \cite{M2020} gave a general upper bound on the star chromatic index of
biregular graphs.
\begin{thm}[\cite{M2020}]
Let $G$ be an $(r,d)$-biregular graph, where $r\geq d>1$. Then $\chi'_{st}(G)\leq r^2-2r+d+1$.
\end{thm}

 Now we  consider the star chromatic index of bipartite
graphs in terms of the maximum degree. By Theorem \ref{s=7}, we have $\chi'_{st}(G)\leq7$ if $G$ is a bipartite graph with $\Delta(G)=3$. Wang, Wang, and Wang \cite{WWW2019} considered bipartite graphs with $\Delta(G)=4$ and proved the following.
\begin{thm}[\cite{WWW2019}]
Let $G$ be a bipartite graph with $\Delta(G)=4$. Then $\chi'_{st}(G)\leq13$.
\end{thm}
Since $K_{4,4}$ requires 7 colours, and no graph requiring 8 colours has been
found, the smallest upper bound is at least 7 \cite{CGR2019}.

Recently, Melinder \cite{M2020} gave a general upper bound for bipartite graphs with maximum degree.
\begin{thm}[\cite{M2020}]
If $G$ is a bipartite graph with maximum degree $\Delta$, then $$\chi'_{st}(G)\leq \Delta^2-\Delta+1.$$
\end{thm}

Melinder \cite{M2020} also proved an upper bound for a special case of multipartite graphs.
\begin{thm}[\cite{M2020}]
Let $K_{1,\ldots,1,n}$ be a complete multipartite graph with $m$ parts of size one. Then
$$\chi'_{st}(K_{1,\ldots,1,n})\leq\begin{cases}
		n      & if~m=1,\\
		2n-\lfloor\frac{n}{2}\rfloor+1   & if~m=2,\\
        3\lceil\frac{n}{2}\rceil+3   & if~m=3,\\
        20\lceil\frac{n}{12}\rceil+6   & if~m=4,\\
        15\lceil\frac{n}{8}\rceil\lceil\frac{m}{8}\rceil+\frac{m(m-1)}{2}   & if~m\geq5.\\
	\end{cases}
$$
\end{thm}

\section{Planar graphs}\label{Planargraphs}
Deng, Liu, and Tian  \cite{DLT2011} and Bezegov\'{a}, Lu{\v{z}}ar, Mockov{\v{c}}iakov{\'a}, Sot{\'a}k,  and
  {\v{S}}krekovski \cite{BLM2016} independently proved that for each tree $T$ with
maximum degree $\Delta$, its star chromatic index $\chi'_{st}(T)\leq\left\lfloor\frac{3\Delta}2\right\rfloor$,  Moreover,  the bound is tight.
Bezegov\'{a}, Lu{\v{z}}ar, Mockov{\v{c}}iakov{\'a}, Sot{\'a}k,  and
  {\v{S}}krekovski \cite{BLM2016} investigated the star edge-coloring of outerplanar graphs by showing the following results.

\begin{thm}[\cite{BLM2016}]\label{outerplanar}
Let $G$ be an outerplanar graph. Then
\begin{enumerate}[(a)]
\item   $\chi'_{st}(G) \le\left\lfloor\frac{3\Delta(G)}{2}\right\rfloor+12$  if $\Delta(G)\ge4$.
\item  $\chi'_{st}(G)\le5$ if $\Delta(G)\le3$.
\end{enumerate}
\end{thm}

In \cite{BLM2016}, the authors pointed out that the constant 12 in Theorem \ref{outerplanar}(a) can be decreased to 9 by using more involved analysis. Moreover, they put forward the following conjecture.
\begin{conj}[\cite{BLM2016}]\label{outerplanar1}
 Every outerplanar graph $G$ has $\chi'_{st}(G) \le\left\lfloor\frac{3\Delta(G)}{2}\right\rfloor+1$  if $\Delta(G)\ge4$.
\end{conj}

Wang, Wang, and Wang \cite{WWW2018} improved  Theorem \ref{outerplanar}(a) by edge partition method.
\begin{thm}[\cite{WWW2018}]\label{outerplanar5}
If $G$ is an outerplanar graph, then $\chi'_{st}(G) \le\left\lfloor\frac{3\Delta(G)}{2}\right\rfloor+5$.
\end{thm}

A {\em cactus} is a graph in which every edge belongs to at most one cycle. Since these graphs are outerplanar, in order
to prove Conjecture \ref{outerplanar1}, it is worth to study the star edge
coloring of cactus graphs. Omoomi, Dastjerdi, and Yektaeian \cite{ODY2020} proved Conjecture
\ref{outerplanar1} for cactus graphs.
\begin{thm}[\cite{ODY2020}]
If $G$ is a cactus, then $\chi'_{st}(G)\leq\left\lfloor\frac{3\Delta(G)}{2}\right\rfloor+1$.
\end{thm}

An outerplanar graph $G$ is {\em maximal} if $G+uv$ is not outerplanar for any two non-adjacent vertices $u$
and $v$ of $G$. Deng and Tian \cite{DT2011} determined the exact values of the star chromatic index for all maximal outerplanar graphs with order $3\leq n\leq8$ and proved the following results.
\begin{thm}[\cite{DT2011}]
Let $G$ be a maximal outerplanar graph with order $n$ and maximum degree $\Delta$. Then
\begin{itemize}
  \item $\chi'_{st}(G)=6$ if $\Delta=4$.
  \item $6\leq\chi'_{st}(G)\leq n-1$ if $n\geq8$. The bounds are tight.
\end{itemize}
\end{thm}

Recently, Deng, Yao, Zhang, and Cui \cite{DYZC2020} studied the star chromatic index of 2-connected outerplanar graphs and proved that if
$G$ is a 2-connected outerplanar graph with diameter $d$ and maximum degree $\Delta$, then
$\chi'_{st}(G) \le\Delta+6$ if $d=2$ or $3$;  $\chi'_{st}(G) \le9$ if $\Delta=5$.


Wang, Wang, and Wang  \cite{WWW2018} investigated the star edge-coloring of planar graphs by using an edge-partition technique and a useful relation between star chromatic index and strong chromatic index.
\begin{defn}
A proper $k$-edge-coloring of $G$ is called a {\em strong $k$-edge-coloring} if any two edges of distance at most two get distinct colors. That is, each color class is an induced matching in the graph $G$. The {\em strong chromatic index}, denoted $\chi'_{s}(G)$, of $G$ is the smallest integer $k$ such that $G$ has a strong $k$-edge-coloring.
\end{defn}

From the definitions, it is straightforward to see that $\chi'_{st}(G)\leq\chi'_{s}(G)$.
\begin{defn}
Suppose that $H$ is a subgraph of a graph $G$. A {\em restricted-strong} $k$-edge-coloring of $H$ on $G$ is a $k$-edge-coloring  such that any two edges  having a distance at most two in $G$ get distinct colors. The {\em restricted strong chromatic index} of $H$ on $G$, denoted $\chi'_{s}(H|_G)$, is the smallest integer $k$ such that $H$ has a restricted strong $k$-edge-coloring on $G$.
\end{defn}
\begin{defn}
An {\em edge-partition} of $G$ is a decomposition of $G$ into subgraphs $G_1, G_2,\ldots,G_m$ such that $E(G)=E(G_1)\cup E(G_2) \cup\cdots E(G_m)$ and $E(G_i)\cap E(G_j)=\emptyset$ for $i\neq j$.
\end{defn}
\begin{thm}[\cite{WWW2018}]\label{keylemma}
If a graph $G$ can be edge-partitioned into two graphs $F$ and $H$, then
$$\chi'_{st}(G)\leq\chi'_{st}(F)+\chi'_{s}(H|_G).$$
\end{thm}
\begin{thm}[\cite{WHW2017,WWW2018}]\label{edgepartition}
Every planar graph $G$ with maximum degree $\Delta$ can be edge-partitioned into two forests $F_1$, $F_2$ and a subgraph $K$ such that $\Delta(K)\leq10$ and $\Delta(F_i)\leq\lceil(\Delta-9)/2\rceil$ for $i=1,2$
\end{thm}

Using Theorems \ref{keylemma} and \ref{edgepartition}, currently the best bound for
planar graphs can be obtained.\begin{thm}[\cite{WWW2018}]
Let $G$ be a planar graph with maximum degree $\Delta$. Then
$$\chi'_{st}(G) \le2.75\Delta+18.$$
\end{thm}
Similarly, Wang, Wang, and Wang  \cite{WWW2018} proved more specific results (together with the
result for outerplanar graphs from Theorem \ref{outerplanar5}).
\begin{thm}[\cite{WWW2018}]
Let $G$ be a planar graph with maximum degree $\Delta$ and girth $g$. Then
\begin{enumerate}[(a)]\label{wangplanar}
\item $\chi'_{st}(G) \le2.25\Delta+6$ if $G$ is $K_4$-minor free.
\item $\chi'_{st}(G) \le\lfloor1.5\Delta\rfloor+18$ if $G$ has no $4$-cycles.
\item $\chi'_{st}(G) \le\lfloor1.5\Delta\rfloor+13$ if $g\geq5$.
\item $\chi'_{st}(G) \le\lfloor1.5\Delta\rfloor+3$ if $g\geq8$.
\end{enumerate}
\end{thm}

\section{List star edge-coloring of graphs}\label{Liststaredgecoloring}
A natural generalization of star edge-coloring is the list star edge-coloring.
Given a list assignment $L$ which assigns to each edge $e$ a finite set $L(e)$, a graph is said to be {\em $L$-star edge-colorable} if $G$ has a star edge-coloring $c$ such that $c(e)\in L(e)$ for each edge $e$. The {\em list-star chromatic index}, $ch'_{st}(G)$ of a graph $G$ is  the minimum $k$ such that for every edge list $L$ for $G$ with $|L(e)|=k$ for every $e\in E(G)$, $G$ is $L$-star edge-colorable.
For any graph $G$, it is obvious that $\chi'_{st}(G)\leq ch'_{st}(G)$.

Moser and Tardos \cite{MT2010} designed an algorithmic version of Lov\'{a}sz Local Lemma by means of the so-called {\em Entropy Compression Method}.
In \cite{CYY2019}, by employing Entropy Compression Method, Cai, Yang, and Yu gave an upper bound for the list-star chromatic index of a graph. \begin{thm}[\cite{CYY2019}]
 For any graph $G$ with maximum degree $\Delta$, $$ch'_{st}(G)\leq\left\lceil2\Delta^{\frac{3}{2}}(\frac{1}{\Delta}+2)^\frac{1}{2}+2\Delta\right\rceil.$$
 \end{thm}

Dvo\v{r}\'ak, Mohar, and \v{S}\'amal~\cite{DMS2013} asked the following question.
\begin{que}[\cite{DMS2013}]\label{list}
Is it true that $ch'_{st}(G)\leq7$ for every
subcubic graph $G$? (perhaps even $\leq6$?)
\end{que}

Kerdjoudj, Kostochka, and Raspaud~\cite{KKP2017} gave a partial answer to this question. They proved that $ch'_{st}(G)\leq8$ for  every subcubic graph. Subsequently, Lu\v{z}ar, Mockov\v{c}iakov\'{a}, and Sot\'{a}k \cite{LMS2019} answered Question \ref{list} in affirmative.
\begin{thm}[\cite{LMS2019}]
Let G be a subcubic graph. Then $ch'_{st}(G)\leq7$.
\end{thm}

Kerdjoudj, Kostochka, and Raspaud \cite{KKP2017}, Kerdjoudj and Kostochka \cite{KP2018} and Kerdjoudj, Pradeep,
and Raspaud \cite{KPP2018} studied the list star edge-coloring of graphs with small maximum average degree.

\begin{thm} [\cite{KKP2017,KP2018,KPP2018}]\label{girth5} Let $G$ be a  graph with maximum degree $\Delta$.
\begin{enumerate}[(a)]
\item If $\text{mad}(G)<7/3$, then $ch'_{st}(G)\leq 2\Delta-1$.
\item If $\text{mad}(G)<5/2$, then $ch'_{st}(G)\leq 2\Delta$.
\item If $\text{mad}(G)<8/3$, then $ch'_{st}(G)\leq 2\Delta+1$.
\item If $\text{mad}(G)<14/5$, then $ch'_{st}(G)\leq 2\Delta+2$.
\item If $\text{mad}(G)<3$, then $ch'_{st}(G)\leq 2\Delta+3$.
\end{enumerate}
\end{thm}

In \cite{LHLK2020}, Li, Horacek, Luo, and  Miao presented, in contrast,
a best possible linear upper bound for the list-star chromatic index for some graph classes.
Specifically they showed the following.
\begin{thm}[\cite{LHLK2020}]\label{girth83}
Let $\varepsilon>0$ be a real number and $d=2\lceil\frac{8-3\varepsilon}{9\varepsilon}\rceil$. Let $G$ be a graph with maximum degree $\Delta$ and  $mad(G)<\frac{8}{3}-\varepsilon$. Then
$$ch'_{st}(G)\leq\max\{\frac{3}{2}\Delta+\frac{d}{2}+2,\Delta+2d+1\}.$$
\end{thm}
To overcome some difficulties in the proof, they developed a new
coloring extension method by requiring certain edges to be colored with different size of sets of
colors.

 It was observed in~\cite{girth} that  every planar graph with girth $g$ satisfies $\text{mad}(G)< \frac {2g}{g-2}$.   This implies the following by Theorems \ref{mainthm}, \ref{girth5}, and \ref{girth83}.

\begin{cor}[\cite{KKP2017,KP2018,KPP2018,LSSW2018,LHLK2020}]\label{listplanar}
Let $G$ be a planar graph with maximum degree $\Delta$ and girth $g$.
\begin{enumerate}[(a)]
\item  If $\Delta=3$ and $g\geq12$, then $\chi'_{st}(G)\leq 5$.
\item If $g\geq14$, then $ch'_{st}(G)\leq 2\Delta-1$.
\item If $g\geq10$, then $ch'_{st}(G)\leq 2\Delta$.
\item If $g\geq8$, then $ch'_{st}(G)\leq 2\Delta+1$.
\item If $g\geq7$, then $ch'_{st}(G)\leq 2\Delta+2$.
\item If $g\geq6$, then $ch'_{st}(G)\leq 2\Delta+3$.
\item If $g\geq9$, then $ch'_{st}(G)\leq \max\{\frac{3\Delta}{2}+11, \Delta+37\}$.
\item If $g\geq16$, then $ch'_{st}(G)\leq \frac{3\Delta}{2}+4$.
\end{enumerate}
\end{cor}

For planar graphs with girth $4$,  Li, Horacek, Luo, and  Miao \cite{LHLK2020} obtained an infinite family of such graphs whose list-star
chromatic index can not be bounded above by $\frac{3}{2}\Delta+c$.\begin{prop}[\cite{LHLK2020}]\label{girth4}
For each integer $\Delta\geq3$, there exists a planar graph $G$ of girth 4 with
maximum degree $\Delta$  such that $$ch'_{st}(G)\geq \chi'_{st}(G)\geq\frac{13}{8}\Delta-\frac{3}{4}.$$
\end{prop}
In view of Theorem \ref{wangplanar}(c), Corollary \ref{listplanar}(g) and Proposition \ref{girth4},  Li, Horacek, Luo, and  Miao \cite{LHLK2020} conjectured that
\begin{conj}[\cite{LHLK2020}]
There exists a constant $c>0$ such that for any planar graph $G$ of girth at
least 5 with maximum degree $\Delta$, we have $$ch'_{st}(G)\leq\frac{3}{2}\Delta+c.$$
\end{conj}

Let $k\in \mathbb{N}$. A graph $G$ is $k$-degenerate if $\delta(H)\leq k$ for every subgraph $H$ of $G$. Kerdjoudj and Raspaud \cite{KP2019} and Han, Li, Luo, and Miao \cite{HLLM2019} proved:
\begin{thm}[\cite{KP2019}]\label{3k-2}
Every $k$-degenerate graph $G$ of maximum
degree $\Delta$ and $k\geq2$ satisfies
$$ch'_{st}(G)\leq (3k-2)\Delta-k^2+2.$$
\end{thm}

\begin{thm}[\cite{HLLM2019}]
 Let $k\geq2$ be an integer. For every $k$-degenerate graph $G$ with maximum
degree $\Delta$ , we have the following two upper bounds:
\begin{enumerate}[(a)]
\item   $ch'_{st}(G)\leq \frac{5k-1}{2}\Delta-\frac{k(k+3)}{2}$. The bound is tight for $C_5$ as $ch'_{st}(C_5)=4$
\item  $ch'_{st}(G)\leq 2k\Delta+k^2-4k+2$.
\end{enumerate}
\end{thm}

Kerdjoudj and Raspaud \cite{KP2019} also considered the class of $K_4$-minor free graphs. Since
the $K_4$-minor free graphs are 2-degenerate\cite{D1965}, by applying Theorem \ref{3k-2} we have $ch'_{st}(G)\leq 4\Delta-2$ for any graph $G$ which is $K_4$-minor free and contains at least one edge.
Kerdjoudj and Raspaud \cite{KP2019} improved this  bound by
proving the following.
\begin{thm}[\cite{KP2019}]\label{3k-3minor}
Let $G$ be a $K_4$-minor free graph with maximum
degree $\Delta$. Then
$$ch'_{st}(G)\leq\begin{cases}
		4      & if~\Delta=2,\\
		3\Delta-3   & if~\Delta\geq3.\\
	\end{cases}
$$
\end{thm}

Notice that Theorem \ref{3k-3minor} implies that Conjecture \ref{cubic} is true for every $K_4$-minor free subcubic graph.

Han, Li, Luo, and Miao \cite{HLLM2019} extended this result of the star chromatic index of trees  to list-star chromatic index.
\begin{thm}[\cite{HLLM2019}]
 For every tree $T$ with maximum degree $\Delta$, $ch'_{st}(T)\leq \lfloor\frac{3\Delta}{2}\rfloor$. This bound is tight.
\end{thm}
\section{Another family of graphs}\label{Anotherniceclasses}
\subsection{Paths, Cycles, Fan graphs, Wheel graphs}
For any integer $n\geq4$, the {\em wheel graph}  $W_n$ (resp. {\em fan graph} $F_n$) is the $n$-vertex graph obtained
by joining a vertex  to each of the $n-1$ vertices
of the cycle graph $C_{n-1}$ (resp. the path graph $P_{n-1}$).
The following results were gave by Deng \cite{D2007} and Wang, Wang, and Wang \cite{WWW2019}.
\begin{thm}[\cite{D2007,WWW2019}]
For $P_n$ ($n\geq2$) and $C_n$ ($n\geq3$),
\begin{enumerate}[(a)]
  \item $ch'_{st}(P_n)=\chi'_{st}(P_n)=\begin{cases}
		1      & if~n=2,\\
		2     & if~3\leq n\leq4,\\
        3    & if~n\geq5.
\end{cases}	
$
  \item $ch'_{st}(C_n)=\chi'_{st}(C_n)=\begin{cases}
		3      & if~n\neq5,\\
        4    & if~n=5.
        \end{cases}
        $
\end{enumerate}

\end{thm}

\begin{thm}[\cite{D2007}]
Let $F_n$ be a fan graph on order $n\geq3$ and $W_n$ be a wheel graph on order $n\geq4$. Then
\begin{enumerate}[(a)]
  \item $\chi'_{st}(F_n)=\begin{cases}
		n+1      & if~n=5,\\
		n     & if~n=3,4,6,7,\\
        n-1    & if~n\geq8.
\end{cases}	
$
  \item $\chi'_{st}(W_n)=\begin{cases}
		n+2      & if~n=5,\\
        n+1    & if~n=4,6,7,\\
        n-1     & if~n\geq8.
        \end{cases}
        $
\end{enumerate}

\end{thm}
\subsection{Halin  graphs}\label{Halingraphs}
A {\em Halin graph} is a planar graph that consists of a plane embedding of a tree $T$ and a cycle $C$, in which the cycle
$C$ connects the leaves of the tree such that it is the boundary of the exterior face and the degree of each
interior vertex of $T$ is at least three.
A {\em complete Halin} graph is a graph that all leaves of the characteristic tree are at the same distance
from the root vertex.

A {\em caterpillar} is a tree whose removal of leaves results in a path called the {\em spine} of the caterpillar. The {\em necklace graph}  denoted by $\mathcal{N}_n$ is a cubic Halin graph  obtained by
joining a cycle with all vertices of degree 1 of a caterpillar  having $n$ vertices of
degree 3 and $n+2$ vertices of degree 1.

  Hou, Li, and Wang \cite{HLW2020} studied the star chromatic index of Halin graphs. Theorem \ref{cubichalin} was also proved by Casselgren, Granholm, and Raspaud \cite{CGR2019}
\begin{thm}[\cite{HLW2020,CGR2019}]\label{cubichalin}
If $G$  is a cubic Halin graph, then $$\chi'_{st}(G)\leq6.$$ Furthermore, the bound is tight.
\end{thm}
\begin{thm}[\cite{HLW2020}]
Let $\mathcal{N}_n$ be a necklace such that $n\geq1$ is odd. Then $$\chi'_{st}(\mathcal{N}_n)\leq5.$$
\end{thm}
\begin{thm}[\cite{HLW2020}]
Let $G$ be a complete Halin graph with maximum degree $\Delta(G)\geq6$. Then $$\chi'_{st}(G)\leq\left\lfloor\frac{3\Delta}{2}\right\rfloor+1.$$
\end{thm}
\subsection{$k$-power graphs}\label{powergraphs}
The {\em $k$-power} of a graph $G$ is the graph $G^k$ whose vertex set is $V(G)$, two distinct vertices being adjacent
in $G^k$ if and only if their distance in $G$ is at most $k$.

Hou, Li, and Wang \cite{HLW2020} studied the star chromatic index of path square graphs and cyclic square graphs.
\begin{thm}[\cite{HLW2020}]
For the graph $P^2_n$ with $n\geq3$,
$$\chi'_{st}(P^2_n)=\begin{cases}
		3      & if~n=3,\\
		4     & if~n=4,\\
        6     & if~n\geq5.
	\end{cases}
$$
\end{thm}
\begin{thm}[\cite{HLW2020}]
For the graph $C^2_n$ with $n\geq3$,
$$\chi'_{st}(C^2_n)\leq\begin{cases}
		8     & \text{if $n$ is odd ($n\neq5$ or 9)},\\
		9     &  \text{if $n$ is even}.
	\end{cases}
$$
\end{thm}
\subsection{Generalized Petersen graphs}\label{GeneralizedPetersengraphs}
Let $n$ and $k$ be positive integers, $n\geq2k+1$ and $n\geq3$. The {\em generalized Petersen
graph} $P(n,k)$, which was introduced in \cite{W1969}, is a cubic graph with $2n$ vertices, denoted
by $\{u_0,u_1,\ldots,u_{n-1},v_0,v_1,\ldots,v_{n-1}\}$, and all edges are of the form $u_iu_{i+1}$,
$u_iv_i$, $v_iv_{i+k}$ for $0\leq i\leq n-1$. In the absence of a special claim, all subscripts of
vertices of $P(n,k)$ are taken modulo $n$ in the following.

By Theorem \ref{s=7}(b), we have $\chi'_{st}(P(n,k))\geq4$. Zhu and Shao \cite{ZS2019} gave a necessary and sufficient condition of $\chi'_{st}(P(n,k))=4$ and  showed that ``almost all'' generalized Petersen graphs have a star 5-edge-coloring. Hou, Li, and Wang \cite{HLW2020} proved that $\chi'_{st}(P(3n,n))=5$ for $n\geq2$, which is a special case of Theorems \ref{pnk4} and \ref{pnk5}.
\begin{thm}[\cite{ZS2019}]\label{pnk4}
$\chi'_{st}(P(n,k))=4$ if and only if $n\equiv 0\pmod4$ and $k\equiv 1\pmod2$.
\end{thm}
\begin{thm}[\cite{ZS2019}]\label{pnk5}
Let $\ell$ be the greatest common divisor of $n$ and $k$. Then $P(n,k)$ has a star
5-edge-coloring when
\begin{itemize}
  \item $\ell\geq3$, with the exception of $\ell=3$, $k\neq\ell$, and $\frac{n}{3}\equiv 1\pmod3$,
  \item $\ell=1$, $n\equiv 0\pmod2$ and $k\equiv 1\pmod2$,
  \item $k=1$ and $n\geq5$,
  \item $k=2$ and $n\equiv0\pmod6$.
\end{itemize}
\end{thm}
Furthermore, Zhu and Shao \cite{ZS2019} found that the generalized Petersen graph
$P(n,k)$ with $n=3$, $k=1$ is the only graph with a star chromatic index of 6
among the investigated graphs. Based on these results, they conjectured that $P(3,1)$
is the unique generalized Petersen graph that admits no star 5-edge-coloring.
\subsection{Cartesian product of graphs}
The {\em Cartesian product} of two graphs $G$ and $H$, denoted by $G\square H$, is a graph with vertex set
$V(G)\times V(H)$, and $(a,x)(b,y)\in E(G\square H)$ if either $ab\in E(G)$ and $x=y$, or $xy\in E(G)$ and
$a=b$. A {\em $d$-dimensional grid} $G_{\ell_1,\ell_2,\ldots,\ell_d}=P_{\ell_1}\square P_{\ell_2}\square\cdots \square P_{\ell_d}$ is the Cartesian product of $d$ paths.
A {\em $d$-dimensional hypercube} $Q_d$ is the Cartesian product of $P_2$ by itself $d$ times. A {\em $d$-dimensional
toroidal grid} $T_{\ell_1,\ell_2,\ldots,\ell_d}=C_{\ell_1}\square C_{\ell_2}\square\cdots \square C_{\ell_d}$ is the Cartesian product of $d$ cycles.

Omoomi and Dastjerdi \cite{OD2018} first gave an upper bound on the star chromatic index of the Cartesian product of
two graphs in terms of their star chromatic indices and their chromatic numbers.
\begin{thm}[\cite{OD2018}]\label{starchromatic}
For any two graphs $G$ and $H$,
$$\chi'_{st}(G\square H)=\chi'_{st}(H\square G)\leq\min\{\chi'_{st}(G)\chi(H)+\chi'_{st}(H),\chi'_{st}(H)\chi(G)+\chi'_{st}(G)\}.$$
Moveover, this bound is tight.
\end{thm}

\begin{thm}[\cite{OD2018}]
For every graph $G$ and a positive integer $n$, the following hold.
\begin{itemize}
  \item If $n\geq2$, then $\chi'_{st}(G\square P_n)\leq \chi'_{st}(G\square C_{2n})\leq \chi'_{st}(G)+2\chi(G)$.
  \item If $n\geq2\chi(G)+1$ is odd, then $\chi'_{st}(G\square C_{n})\leq \chi'_{st}(G)+2\chi(G)+1$.
  \item If $n\geq3$ is odd, then $\chi'_{st}(G\square C_{n})\leq \chi'_{st}(G)+2\chi(G)+\lceil\frac{2\chi(G)}{n-1}\rceil\leq \chi'_{st}(G)+2\chi(G)+3$.
\end{itemize}
\end{thm}
Then they determined the exact value of the star chromatic index of 2-dimensional grids and extended this
result to get an upper bound on the star chromatic index of $d$-dimensional grids, where $d\geq3$. The following theorem \ref{pmpncartesian}, corollary \ref{corgrid} and corollary \ref{corhypercube} also were proved by Deng, Liu, and Tian \cite{DLT2012}.
\begin{thm}[\cite{DLT2012,OD2018}]\label{pmpncartesian}
For all integers $2\leq m\leq n$,
$$\chi'_{st}(P_m\square P_n)=\begin{cases}
		3      & if~m=n=2,\\
		4     & if~m=2,n\geq3,\\
        5     & if~m=3, n\in\{3,4\},\\
        6     & otherwise.
	\end{cases}
$$
\end{thm}

\begin{cor}[\cite{DLT2012,OD2018}]\label{corgrid}
If $G_{\ell_1,\ell_2,\ldots,\ell_d}$ is a  $d$-dimensional grid, where $d\geq2$, then $$\left\lceil\frac{8}{3}(d-\sum\limits_{i=1}^d\frac{1}{\ell_i})\right\rceil\leq\chi'_{st}(G_{\ell_1,\ell_2,\ldots,\ell_d})\leq 4d-2.$$  Moreover, for $d=2$ and $\ell_1,\ell_2\geq4$, the bound is tight.
\end{cor}

In the following theorem, Omoomi and Dastjerdi \cite{OD2018} considered the star chromatic index of the Cartesian product of paths and cycles.
\begin{thm}[\cite{OD2018}]\label{pathscycles}
For all integers  $m\geq2$ and $n\geq3$,
$$\chi'_{st}(P_m\square C_n)\leq\begin{cases}
		7      & \text {if $n$ is even},\\
		8     & \text {if $n$ is odd}.
	\end{cases}
$$
\end{thm}

Now, by Theorems \ref{starchromatic} and \ref{pathscycles}, Omoomi and Dastjerdi \cite{OD2018} obtained an upper bound on the star chromatic index of hypercube $Q_d$, where $d\geq3$, as follows.
\begin{cor}[\cite{DLT2012,OD2018}]\label{corhypercube}
If $Q_d$ is a  $d$-dimensional hypercube, where $d\geq3$, then $$\left\lceil\frac{4}{3}d\right\rceil\leq\chi'_{st}(Q_d)\leq 2d-2.$$ Moreover the bound is tight.
\end{cor}

As a corollary of Theorem \ref{pmpncartesian}, Holub, Lu\v{z}ar, Mihalikov\'{a},  Mockov\v{c}iakov\'{a}, and Sot\'{a}k \cite{HLMMS2020} established the lower bound of 6 colors for the Cartesian products,  where one factor is a cycle and the other is a path of length at least 2. Then they also obtained  the exact values of the star chromatic index for specific lengths of paths and cycles.
\begin{cor}[\cite{HLMMS2020}]\label{pmcn}
For all integers  $m,n\geq3$,
$$\chi'_{st}(P_m\square C_n)\geq6.$$
\end{cor}
\begin{thm}[\cite{HLMMS2020}]
For all integers  $m\geq2$ and $n\geq3$,
$$\chi'_{st}(P_m\square C_n)=\begin{cases}
        4          & \text {if $m=2$ and $n\equiv0\pmod4$},\\
        5       & \text {if $m=2$ and $n\neq3$, $n\not\equiv 0\pmod4$},\\
		6      & \text {if $m=2$, $n=3$ or $m\in\{3,4,5,6\}$ or $m\geq7$, $n\equiv0\pmod{2,3}$ },\\
		7     & \text {if $m\geq7$ and $n\in\{5,7\}$}.
	\end{cases}
$$
\end{thm}

It remains to determine  $\chi'_{st}(P_m\square C_n)$ for $m\geq7$ and $n\geq11$, $n\not\equiv0\pmod{2,3}$. Holub, Lu\v{z}ar, Mihalikov\'{a},  Mockov\v{c}iakov\'{a}, and Sot\'{a}k \cite{HLMMS2020} proposed the following conjecture.
\begin{conj}[\cite{HLMMS2020}]\label{pmcn6}
There exist constants $c_1$ and $c_2$ such that for all integers $m\geq c_1$ and $n\geq c_2$, we have
$$\chi'_{st}(P_m\square C_n)=6.$$
\end{conj}

Finally, Omoomi and Dastjerdi \cite{OD2018} gave some upper bounds on the star chromatic
index of the Cartesian product of two cycles and $d$-dimensional toroidal grids as follows.
\begin{thm}[\cite{OD2018}]
For all integers  $m,n\geq3$,
$$\chi'_{st}(C_m\square C_n)\leq\begin{cases}
		7      & \text {if $m$ and $n$ are even},\\
		8     & \text {if $m\neq3$ is odd and $n$ is even},\\
        9     & \text {if $m=3$ and $n$ is even},\\
        10     & \text {if $m$ and $n$ are odd}.
	\end{cases}
$$
\end{thm}

\begin{cor}[\cite{OD2018}]
For all integers $d\geq2$ and $\ell_1,\ell_2,\ldots,\ell_d\geq3$, $$\chi'_{st}(T_{2\ell_1,2\ell_2,\ldots,2\ell_d})\leq 4d-1~~\text{and}~~\chi'_{st}(T_{\ell_1,\ell_2,\ldots,\ell_d})\leq 7d-4.$$
\end{cor}

Holub, Lu\v{z}ar, Mihalikov\'{a},  Mockov\v{c}iakov\'{a}, and Sot\'{a}k \cite{HLMMS2020}  obtained  the exact values of the star chromatic index for specific lengths of  cycles and showed the star chromatic index of the Cartesian product of two cycles is at most 7. Note that Corollary \ref{pmcn} implies that the Cartesian product of any two cycles will need at least 6 colors for a star edge-coloring.
\begin{thm}[\cite{HLMMS2020}]
For all integers $m,n\geq3$,
$$6\leq\chi'_{st}(C_m\square C_n)\leq7.$$
\end{thm}
\begin{thm}[\cite{HLMMS2020}]
For all integers $m,n\geq3$,
$$\chi'_{st}(C_m\square C_n)=\begin{cases}
		6      & \text {if $m,n\equiv 0\pmod3$ or $m\equiv 0\pmod4$, $n\equiv 0\pmod 2$} \\&\text { or $m=6$, $n\equiv0\pmod{3,4}$},\\
		7     & \text {if $m=3$, $n\not\equiv0\pmod3$ or $m=4$, $n\in\{5,7,9,11\}$} \\& \text { or $m\in\{5,7\}$ or $m=6$, $n\not\equiv0\pmod{3,4}$ or $m=8$, $n=9$}.
	\end{cases}
$$
\end{thm}
Note that Conjecture \ref{pmcn6} is equivalent to the following.
\begin{conj}[\cite{HLMMS2020}]\label{cmcn6}
There exists a constant $c$ such that for every integer $m\geq c$, there exists an integer $n$ such that
$$\chi'_{st}(C_m\square C_n)=6.$$
\end{conj}
In fact, Holub, Lu\v{z}ar, Mihalikov\'{a},  Mockov\v{c}iakov\'{a}, and Sot\'{a}k \cite{HLMMS2020} believed (with a bit lower confidence) that the following stronger version of Conjecture \ref{cmcn6} can also be confirmed.
\begin{conj}[\cite{HLMMS2020}]
There exists a constant $c$ such that for all integers  $m,n\geq c$,
$$\chi'_{st}(C_m\square C_n)=6.$$
\end{conj}
\subsection{Corona product of graphs}\label{Coronaproduct}
The {\em corona product} of two graphs $G$ and $H$ is the graph $G\circ H$ formed from one copy of
$G$ and $|V(G)|$ copies of $H$ where the $i^{th}$ vertex of $G$ is adjacent to every vertex in the
$i^{th}$ copy of $H$. We first state some special graph classes as follows.
\begin{itemize}
  \item The {\em $n$-sunlet} graph on $2n$ vertices is obtained by attaching $n$ pendant edges to the cycle
$C_n$ and is denoted by $C'_n$.
  \item {\em Double star} $K_{1,n,n}$ is a tree obtained from the star $K_{1,n}$ by adding a new pendant edge of the existing $n$ pendant vertices. It has $2n+1$ vertices and $2n$ edges.

  \item The {\em helm graph} $H_n$ is the graph obtained from the wheel graph $W_{n+1}$
by adjoining a pendent edge at each vertices of the outer cycle.
\item The {\em gear graph}  $G_n$, also known as a bipartite wheel graph, is the wheel graph $W_{n+1}$ with a  vertex added between each pair of adjacent vertices of the outer cycle.
\end{itemize}

In \cite{KSV2016} and \cite{KSV2018}, the authors obtained
the star edge chromatic number of the corona product of paths with
paths, paths with
cycles, paths with wheels, paths with helms, paths with gear graphs, paths with
$n$-sunlet graphs, paths with Double star graphs, and paths with complete bipartite graphs,
respectively.

\begin{thm}[\cite{KSV2016}]
For all positive integers $m$ and $n$,
$$\chi'_{st}(P_m\circ P_n)=\begin{cases}
		n      & if~m=1,\\
		n+1     & if~m=2,\\
        n+2     & if~m\geq3.
	\end{cases}
$$
\end{thm}
\begin{thm}[\cite{KSV2018}]
Let $G\in\{P_m\circ C_n, P_m\circ W_n,P_m\circ H_n,P_m\circ G_n\}$, where $m\geq1$ and $n>4$. Then
$$\chi'_{st}(G)=\Delta(G).$$
\end{thm}

\begin{thm}[\cite{KSV2016}]
For all positive integers $m$ and $n\geq3$, then
$$\chi'_{st}(P_m\circ C'_n)=\begin{cases}
		2n      & if~m=1,\\
		2n+1     & if~m=2,\\
        2n+2     & if~m\geq3.
	\end{cases}
$$
\end{thm}
\begin{thm}[\cite{KSV2016}]
For all positive integers $m$ and $n$, then
$$\chi'_{st}(P_m\circ K_{1,n,n})=\begin{cases}
		2n+1      & if~m=1,\\
		2n+2     & if~m=2,\\
        2n+3     & if~m\geq3.
	\end{cases}
$$
\end{thm}
\begin{thm}[\cite{KSV2016}]
For all positive integers $\ell\geq3$, $m\geq3$ and $n\geq3$, then
$$\chi'_{st}(P_\ell\circ K_{m,n})=m+n+2.$$
\end{thm}

\subsection{Graph join of two graphs}
Let $G$ and $H$ be two graphs with $V(G)\cap V(H)=\emptyset$. The {\em graph join} $G\vee H$ of $G$ and $H$ has $V(G\vee H)=V(G)\cup V(H)$ and two different vertices $u$ and $v$ of $G\vee H$ are adjacent if $u\in V(G)$ and $v\in V(H)$, or $uv\in E(G)$ or $uv\in E(H)$. Yang, Liu, and Chen \cite{YLC2008} considered the star chromatic index of the graph join of paths and paths and proved the following results.
\begin{thm}[\cite{YLC2008}]
For any two graphs $G$ and $H$, $$\chi'_{st}(G\vee H)\leq \max\{\chi'_{st}(G),\chi'_{st}(H)\}+|G||H|.$$
\end{thm}
\begin{thm}[\cite{YLC2008}]
For any integer $n\geq2$,
$$\chi'_{st}(P_n\vee P_2)=\begin{cases}
		5      & if~n=2,\\
		\frac{3n-1}{2}+4     & if~n\geq3~and~n\equiv1\pmod2,\\
        \frac{3n}{2}+3     & if~n\geq4~and~n\equiv0\pmod2.
	\end{cases}
$$
\end{thm}
\begin{thm}[\cite{YLC2008}]
$\chi'_{st}(P_3\vee P_3)=10$ and $\chi'_{st}(P_n\vee P_n)<n^2-3n+9$ when $n\geq4$.
\end{thm}
\begin{cor}[\cite{YLC2008}]
$m+2\leq\chi'_{st}(P_m\vee P_n)<mn-3n+9$ when $m>n\geq4$.
\end{cor}

  \vspace{5mm} \noindent {\bf Acknowledgments.}  Hui Lei was partially supported by the National
Natural Science Foundation of China and the Fundamental Research Funds
for the Central Universities, Nankai University.  Yongtang Shi was partially supported by the National
Natural Science Foundation of China (No. 11922112), the Natural Science Foundation of Tianjin and the Fundamental Research Funds
for the Central Universities, Nankai University. 

\end{document}